\documentclass[12pt]{amsart}
\usepackage{amsmath,amsthm,latexsym,amscd,amsbsy,amssymb,amsfonts,
euscript}

\renewcommand{\AE}{\mathop{\rm AE}\nolimits}

\def\ed{\text{e}-\dim}

\newcommand{\Z}{\mathbb Z}

\newtheorem{thm}{Theorem}[section]

\newtheorem{problem}{Problem}

\theoremstyle{definition}
\newtheorem{dfn}{Definition}[section]
\theoremstyle{remark}

\chardef\bslash=`\\ 

\makeatletter
\def\verbatim{\interlinepenalty\@M \@verbatim
  \leftskip\@totalleftmargin\advance\leftskip2pc
  \frenchspacing\@vobeyspaces \@xverbatim}
\makeatother
\numberwithin{equation}{section}

\begin{document}


\title[On two problems in extension theory]
{On two problems in extension theory}
\author{A.~V.~Karasev}
\address{Computer Science and Mathematics Department, Nipissing
University,
100 College Drive, Box 5002, North Bay, Ontario, P1B 8L7, Canada}

\email{alexandk@nipissingu.ca}

\keywords{Extension dimension, universal spaces, quasi-finite
complexes}

\subjclass{Primary 55M10; Secondary 54F45}


\begin{abstract}
In this note we introduce the concept of a quasi-finite complex.
Next, we show that for a given countable and locally finite CW
complex $L$ the following conditions are equivalent:

\begin{itemize}

\item[-] $L$ is quasi-finite.

\item[-] There exists a $[L]$-invertible mapping of a metrizable
compactum $X$ with $\ed X \le [L]$ onto the Hilbert cube.

\end{itemize}

\noindent Finally, we construct an example of a quasi-finite
complex $L$ such that its extension type $[L]$ does not contain a
finitely dominated complex.

\end{abstract}

\maketitle \markboth{A.~V.~Karasev}{On two problems in extension
theory}


\section{Introduction}

 In \cite{ch99} Chigogidze stated the following
two problems and showed that they are equivalent:

\begin{problem}
{\rm Characterize connected locally compact simplicial complexes
$P$ such that $P\in AE (X)$ iff $P\in AE(\beta X)$ for any space
$X$.}
\end{problem}

\begin{problem}
{\rm Characterize connected locally compact simplicial complexes
$P$ such that there exists a $P$-invertible map $f\colon X\to
I^{\omega}$ where $X$ is a metrizable compactum with $P\in AE
(X)$.}
\end{problem}

The following problem in extension theory is closely related to
the Problems 1 and 2 \cite{ch02,dydak96}:

\begin{problem}\label{univ}
{\rm Let $L$ be a countable CW complex such that the class of
metrizable compacta $\{X\colon L\in AE (X)\}$ has a universal
space. Is it true that the extension type $[L]$ of this complex
contains a finitely dominated complex?}
\end{problem}

In this note we introduce the notion of quasi-finite complexes and
show that the class of quasi-finite complexes yields the
characterization required in Problems 1 and 2. Next, we construct
an example of a quasi-finite complex $L$ such that its extension
type $[L]$ does not contain a finitely dominated complex. This
provides a negative solution for Problem 3.

\section{Preliminaries}

For spaces $X$ and $L$, the notation $L \in \AE (X)$ means that
every map $f\colon A \to L$, defined on a closed subspace $A$ of
$X$, admits an extension $\bar f$ over $X$. Let $L$ and $K$ be CW
complexes. Following Dranishnikov \cite{D}, we say that $L \leq K$
if for each space $X$ the condition $L \in {\rm AE}(X)$ implies
the condition $K \in {\rm AE}(X)$. This definition leads to a
preorder relation $\leq$ on the class of CW complexes. This
preorder relation generates the equivalence relation. The
equivalence class of complex $L$ is called the extension type of
$L$ and is denoted by $[L]$. By $\ed X$ we denote extension
dimension of space $X$ \cite{D,DD}. Inequality $\ed X\le [L]$
means that $L\in AE(X)$. More information about extension
dimension and extension types can be found in \cite{ch01}.

The following theorem \cite[Theorem~2.1]{ch99} shows that Problems
1 and 2, stated in the Introduction, are equivalent.

\begin{thm}\label{equivalence} {\rm (A.~Chigogidze).} Let $P$ be a Polish ANR-space. Then
the following statements are equivalent:

\begin{itemize}
    \item[(a)] $P\in AE(\beta X)$ whenever $X$ is a space with $P\in AE (X)$.
    \item[(b)] $P\in AE(\beta X)$ whenever $X$ is a normal space with $P\in AE (X)$.
    \item[(c)] $P\in AE(\beta (\oplus\{ X_t\colon t\in T\} ))$ whenever $T$ is an arbitrary indexing set
    and $X_t$, $t\in T$, is a separable metrizable space with $P\in AE (X_t)$.
    \item[(d)] $P\in AE(\beta (\oplus\{ X_t\colon t\in T\} ))$ whenever $T$ is an arbitrary indexing set
    and $X_t$, $t\in T$, is a Polish space with $P\in AE (X_t)$.
    \item[(e)] There exists a $P$-invertible map $f\colon X\to
I^{\omega}$ where $X$ is a metrizable compactum with $P\in AE
(X)$.
\end{itemize}

\end{thm}

\section{Quasi-finite complexes}

A pair of spaces $V\subset U$ is called {\it $[L]$-connected for
Polish spaces} \cite{brchka}  if for every Polish space $X$ with
$\ed X\le [L]$ and for every closed subspace $A\subset X$ any
mapping of $A$ to $V$ can be extended to a mapping of $X$ into
$U$.

\begin{dfn}
We say that a CW complex $L$ is {\it quasi-finite} if for every
finite subcomplex $P$ of $L$ there exists finite subcomplex $P'$
of $L$ containing $P$ such that the pair $P\subset P'$ is
$[L]$-connected for Polish spaces.
\end{dfn}

\begin{thm}\label{main}
Let $L$ be a countable and locally finite CW complex. Then the
following conditions are equivalent:

\begin{itemize}
    \item[(i)] $L$ is quasi-finite.
    \item[(ii)] Any of the conditions (a)-(e) of Theorem
    \ref{equivalence} is satisfied.
\end{itemize}

\end{thm}

\begin{proof}
It is enough to show that condition (i) is equivalent to the
condition (d) of the Theorem \ref{equivalence}.

Suppose that $L$ is quasi-finite. Let $\{ X_t\colon t\in T\}$ be
an arbitrary family of Polish spaces with $L\in AE (X_t)$ for all
$t\in T$. Consider a closed subspace $A\subset \beta (\oplus\{
X_t\colon t\in T\} )$ and a mapping $f\colon A \to L$. Let $P$ be
a finite subcomplex of $L$ containing $f(A)$.  Since $L$ is
quasi-finite there exists a finite subcomplex $P'$ of $L$
containing $P$ such that the pair $P\subset P'$ is $[L]$-connected
for Polish spaces. Let $\widetilde{f}\colon \widetilde{A} \to P$
be an extension of $f$ over some closed neighborhood
$\widetilde{A}$ of $A$ in $\beta (\oplus\{ X_t\colon t\in T\} )$.
For any $t\in T$, let $f_t \colon X_t\to P'$ be an extension of
$\widetilde{f}|_{\widetilde{A}\cap X_t}$. Consider a mapping
$f'=\oplus \{f_t \colon t\in T\} \colon \oplus\{ X_t\colon t\in
T\}\to P'$. Since $P'$ is compact, mapping $f'$ can be extended to
a mapping $\overline{f} \colon \beta (\oplus\{ X_t\colon t\in T\}
)\to P' \subset L$. Clearly $\overline{f}$ provides a necessary
extension of $f$. Thus $\ed \beta (\oplus\{ X_t\colon t\in T\}
)\le [L]$.

Now suppose that condition (d) of the Theorem \ref{equivalence} is
satisfied. Consider a finite subcomplex $P \subset L$. Let $\{
X_t, A_t, f_t \colon t\in T\}$ be the set of all triples such that
$X_t$ is a Polish space with $\ed X_t \le [L]$, $A_t$ is a closed
subspace of $X_t$ and $f_t\colon A_t \to P$ is a mapping. Put
$f=\oplus \{f_t \colon t\in T\} \colon \oplus\{ A_t\colon t\in
T\}\to P$. Since $P$ is compact there exists a mapping
$\widetilde{f}\colon \overline{\oplus\{ A_t\colon t\in T\}}
^{\,\beta (\oplus\{ X_t\colon t\in T\} )} = \beta (\oplus\{
A_t\colon t\in T\} )\to P$ extending $f$. By assumption $\ed \beta
(\oplus\{ X_t\colon t\in T\} )\le [L]$ and we can extend the
mapping $\widetilde{f}$ to a mapping $\overline{f}\colon \beta
(\oplus\{ X_t\colon t\in T\} )\to L$. Let $P'$ be a finite
subcomplex of $L$ containing $\overline{f} (\beta (\oplus\{
X_t\colon t\in T\} ))$. It is easy to see that the pair $P\subset
P'$ is $[L]$-connected for Polish spaces.
\end{proof}

\section{Example}

In this section we show that there exists a quasi-finite complex
$M$ which is not finitely dominated. By Theorem \ref{main} for
such complex $M$ there exists $[M]$-invertible mapping of a
metrizable compactum $X$ with $\ed X\le [M]$ onto the Hilbert
cube. This implies that $X$ is universal for the class of
metrizable compacta $\{X\colon M\in AE (X)\}$. Thus the solution
of Problem \ref{univ} is negative.

Let $M$ be a countable and locally finite CW complex homotopically
equivalent to the bouquet $\vee \{ M_p \colon p\,\mbox{-prime, }
p\ne 2\}\vee S^3$ where $M_p = M(\Z _p , 2)$ is a Moore space of
the type $(\Z _p,2)$. Clearly $M$ is quasi-finite and $[S^2]\le
[M] < [S^3]$ (further consideration will imply that $[S^2]<[M])$.
We shall show that the extension type $[M]$ of $M$ does not
contain a finitely dominated complex.

Suppose the opposite. Let $L$ be a countable finitely dominated CW
complex such that $[M]=[L]$. Then $[S^2]\le [L] < [S^3]$. In
particular, $L$ is simply connected and therefore $H_2 (L) \cong
\pi _2 (L) $. Since $L$ is finitely dominated, groups $H_2 (L)
\cong \pi _2 (L) $ are finitely generated. Since $[L] < [S^3]$ it
follows that $\pi _2 (L)$ is non-trivial. Consider two cases

\underline{Case 1}. $H_2 (L) \cong\pi _2 (L) \cong\Z\oplus H$
where $H$ is finitely generated.

Note that groups $H^1 (M;\Z _2 )$ and $H^2 (M;\Z _2)$ are trivial.
On the other hand, suspension isomorphism and Hurewicz theorem
imply that $H_3 (\Sigma L ) \cong\Z\oplus H$ and hence the group
$H^ 3(\Sigma L; \Z _2)$ is non-trivial. Therefore we can apply
construction of Dranishnikov and Repov\v s \cite{DR,DR1} and use
the idea from the proof of Theorem 1.4  \cite[p.351]{DR1} to
obtain a metrizable compactum $X$ with $\ed X\le [M]$ and a
mapping $f\colon X\to\Sigma L$ which is not null-homotopic. Let $A
= f^{-1} (L)$ where $L$ is considered as the equator of $\Sigma
L$. Then $f|_A$ does not have an extension $\overline{f}\colon
X\to L$. Indeed, such an extension $\overline{f}$ would be
null-homotopic. On the other hand it would be homotopic to $f$.
This shows that $L\notin AE(X)$ which leads to a contradiction.

\underline{Case 2}. $H_2 (L) \cong\oplus\{\Z _{n_i}\colon
i=1,\cdots ,m\}$

Choose a prime $p> \max\{ n_i \colon i=1,\cdots ,m\}$. Then $H^1
(L , \Z _p)$ and $H^2 (L, \Z _p)$ are trivial. Note that $H^3
(\Sigma M_p , \Z _p)$ is non-trivial. As before we can find a
metrizable compactum $X$ with $\ed X\le [L]$ and a mapping
$f\colon X\to\Sigma M _p$ which is not null-homotopic. By the same
arguments as above we conclude that $M_p\notin AE (X)$. Since
$[M]\le [M_p]$ we obtain a contradiction with $[M]$=$[L]$.

The above consideration shows that extension type $[M]$ does not
contain a finitely dominated complex.


\end{document}